\documentclass[11pt]{article}

\usepackage{amsmath,amssymb,latexsym}
\usepackage{theorem}
\usepackage[dvips]{graphicx}
\usepackage{color}
\newtheorem{theorem}{Theorem}
\newtheorem{lemma}{Lemma}

{\theorembodyfont{\rmfamily} \newtheorem{definition}{Definition}

\newtheorem{remark}{Remark}}
{\theorembodyfont{\slshape} }

\newcommand{\supp}{\mathop{\rm supp}}
\newcommand{\field}[1]{\mathbb{#1}}
\newcommand{\R}{\field{R}}

\newcommand{\N}{\field{N}}
\newcommand{\C}{\field{C}}

\newcommand{\const}{{\rm const}}

\newcommand{\Int}{\mathop{\rm Int}}

\renewcommand{\Re}{\mathop{\rm Re}}

\newcommand{\dist}{\mathop{\rm dist}}

\newcommand{\dsty}{\displaystyle}

{\rm \trivlist \item[\hskip \labelsep{\bf Proof. }]}%
{\hspace*{\fill}$\Box$\endtrivlist}
{\rm \trivlist \item[\hskip \labelsep{\bf Proof}]}%
{\hspace*{\fill}$\Box$\endtrivlist}
\numberwithin{equation}{section}

\title{The Szeg\"{o} Curve and Laguerre polynomials with large
negative parameters}%

\author{C. D\'{\i}az Mendoza\\
Univ. de La Laguna, Spain \\ cjdiaz@ull.es \\
\and R. Orive \\
Univ. de La Laguna, Spain \\ rorive@ull.es}

\date{\today}

\begin{document}

\pagestyle{myheadings} \thispagestyle{plain} \markboth{R.
ORIVE}{SZEG\"{O} CURVE AND LAGUERRE POLYNOMIALS} \maketitle

\begin{abstract}
We study the asymptotic zero distribution of the rescaled Laguerre
 polynomials, $\displaystyle L_n^{(\alpha_n)}(nz)$, with the parameter $\alpha_n$
varying in such a way that $\displaystyle \lim_{n\rightarrow
\infty}\alpha_n/n=-1$. The connection with the so-called Szeg\"{o}
curve will be showed.
\end{abstract}

\section{Introduction}
The definition and many properties of the Laguerre polynomials
$L_n^{(\alpha)}$ can be found in Ch.\ V of Szeg\H{o}'s classic
memoir \cite{szego:1975}. Given explicitly by
\begin{equation}\label{explLag}
L_n^{(\alpha)} (z)=\sum_{k=0}^n
\binom{n+\alpha}{n-k}\frac{(-z)^{k}}{k!}\,,
\end{equation}
or, equivalently, by the well-known Rodrigues formula
\begin{equation}\label{RodrLag}
L_n^{(\alpha)} (z)=\frac{(-1)^n}{n!}\, z^{-\alpha} e^{z}\left(
\frac{d}{dz} \right)^n \left[ z^{n+\alpha} e^{-z}\right]\,,
\end{equation}
they can be considered for arbitrary values of the parameter
$\alpha \in \C$. In particular, (\ref{explLag}) shows that each
$L_n^{(\alpha)}$ depends analytically on $\alpha$ and no degree
reduction occurs: $\deg L_n^{(\alpha)} = n$ for all $\alpha \in
\C$.

For $\alpha>-1$ it is well-known the orthogonality of
$L_n^{(\alpha)}(x)$ on $[0,+\infty)$ with respect to the weight
function $ x^{\alpha} e^{-x}$; in particular, all their zeros are
simple and belong to $[0,+\infty)$. In the general case, $\alpha
\in \C$, $L_n^{(\alpha)}(z)$ may have complex zeros; the only
multiple zero can appear at $z=0$, which occurs if and only if
$\alpha \in \{-1, -2, \dots, -n\}$. In this case we have
\begin{equation}\label{expneg} L_{n}^{(-k)}(z) =
(-z)^{k}\frac{(n-k)!}{n!}\, L_{n-k}^{(k)}(z)\,,
\end{equation} which shows that
$z=0$ is a zero of multiplicity $k$ for $L_{n}^{(-k)}(z)$.

In a series of papers (\cite{Kuijlaars01}, \cite{Kuijlaars04} and
\cite{MR1858305}), asymptotics for rescaled Laguerre polynomials
$\displaystyle L_n^{(\alpha_n)}(nz)$ were analyzed, under the
assumption that $\displaystyle \lim_{n\rightarrow
\infty}\alpha_n/n=A\in \R$. In \cite{MR1858305} the
authors obtained the weak zero asymptotics for the case where
$A<-1$, by means of classical (logarithmic) potential theory. To
this end, it played a key role a full set of non-hermitian
orthogonality relations satisfied by Laguerre polynomials in a
class of open contours in $\C$. Unfortunately, this
analysis could not be extended to the other cases, since for this
approach it is essential the connectedness of the complement to the
support of the asymptotic distribution of zeros (see e.g.
\cite{Gonchar:89} and \cite{Stahl:86}). However, the authors
formulated in \cite{MR1858305} a conjecture for the case $-1<A<0$,
which was proved in some cases and refused in others in \cite
{Kuijlaars04}, by means of the Riemann-Hilbert approach (which has
been previously used by the same authors in \cite{Kuijlaars01} to
obtain strong asymptotics in the case $A<-1$). A similar study for
Jacobi polynomials with varying nonstandard parameters has been
carried out in \cite{ArnoAndrei}, \cite{MR2002d:33017} and
\cite{MFOr05}.

Jacobi or Laguerre polynomials with real parameters (and in
general, depending on the degree n) appear naturally as polynomial
solutions of hypergeometric differential equations, or in the
expressions of the wave functions of many classical systems in
quantum mechanics (see e.g. \cite{Bagrov}).

\begin{figure}[htb] \label{fig:szego}
\centering \includegraphics[scale=0.7]{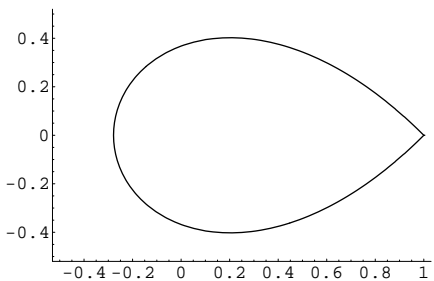} \caption{The
Szeg\H{o} curve.}
\end{figure}

In \cite{MR1858305}, the authors also formulated a conjecture for
the case $A=-1$, but up to now this problem has remained open. Observe
that, by \eqref{expneg}, when $k=n$ we have:
\begin{equation*}\label{multiplezero}
L_{n}^{(-n)}(z)\,=\,(-1)^n\,\frac{1}{n!}\,z^n\,.
\end{equation*}
There is another particular situation corresponding to the case
$A=-1$ which is very well-known in the literature: when
$\alpha_n=-n-1$, we have:
\begin{equation*}\label{exp}
L_{n}^{(-n-1)}(z)\,=\,(-1)^n\,\sum_{k=0}^{n}\frac{z^k}{k!}\,,
\end{equation*}
and thus, in this case the Laguerre polynomials agree (up to a
possible sign) with the partial sums of the exponential series. In
a seminal paper, G. Szeg\H{o} \cite{szego:1924} showed that the
zeros of the rescaled partial sums of the exponential series,
$\displaystyle
\sum_{k=0}^{n}\frac{(nz)^k}{k!}\,=\,(-1)^n\,L_{n}^{(-n-1)}(nz)$,
approach the so-called the Szeg\H{o} curve:
\begin{equation}\label{szego}
\Gamma \,=\,\left\{z\in
\C\,\\,\left|ze^{1-z}\right|=1\,,\,|z|\leq 1\right\}\,,
\end{equation}
\noindent which is a closed curve around the origin passing
through $z=1$ and crossing once the negative real semiaxis
$(-\infty,0)$ (see Figure \ref{fig:szego}). See also \cite{PrVa97}
for a detailed study of the Szeg\"{o} curve and some related problems
in approximation of functions. Recently, T. Kriecherbauer et al.
\cite{KKMM08} obtained uniform asymptotic expansions for the
partial sums of the exponential series by means of the
Riemann-Hilbert analysis. Also, in \cite{DaKu09}, the authors
studied the asymptotics of orthogonal polynomials with respect to
modified Laguerre weights of the type
$$z^{-n+\nu}\,e^{-Nz}\,(z-1)^{2b}\,,$$ \noindent where $n,N\rightarrow
\infty\,$ with $N/n\rightarrow 1\,$ and $\nu\,$ is a fixed number
in $\R\setminus\N\,.$

In this paper, the weak zero asymptotics of rescaled Laguerre
polynomials $\displaystyle L_n^{(\alpha_n)}(nz)$, with
$\displaystyle \lim_{n\rightarrow \infty}\alpha_n/n=-1\,$ will be
analyzed. For it, we will prove that such rescaled Laguerre
polynomials are asymptotically extremal on certain well defined
curves in the complex plane.

The outline of the paper is as follows. In sect 2, the main result
about the weak zero asymptotics of the rescaled Laguerre
polynomials is announced, and in sect. 3, some basic facts on
potential theory and asymptotically extremal polynomials are
recalled. Finally, the proofs are given in sect. 4.

\section{Main Result}

Along with the Szeg\"{o} curve \eqref{szego}, we need to introduce the
family of level curves:
\begin{equation}\label{szegolevel}
\Gamma_r\,=\,\left\{z\in
\mathbb{C}\,\\,\left|ze^{1-z}\right|=e^{-r}\,,\,|z|\leq
1\right\}\,,\,0\leq r< +\infty \,,
\end{equation}
\noindent while for $r=\infty$ we take $\Gamma_{\infty}=\{0\}$.
Observe that $\Gamma_0=\Gamma$, the Szeg\H{o} curve.
We consider the usual counterclockwise orientation. All the level
curves $\Gamma_r$ ($0\leq r<+\infty$) are closed contours such
that $\{0\}\subset \Int(\Gamma_{r})$ and $\Gamma_{r'} \subset
\Int(\Gamma_r)$, for $r'>r$. On the sequel, the interior of
$\Gamma_r$ will be denoted by $G_r$. Associated with this family
of curves, consider for $0\leq r<+\infty$ the family of measures:
\begin{equation}\label{measures}
d\mu_r (z)\, =\,\frac{1}{2\pi i} \frac{1-z}{z} \, dz\,,\,z\in
\Gamma_r\,,
\end{equation}
and set $d\mu_{\infty}(z)=\delta_0$.

Let us recall the definition of balayage (or sweeping out) of a
measure (see e.g. \cite{Saff:97}). Given an open set $\Omega$ with
compact boundary $\partial \Omega$ and a positive measure $\sigma$
with compact support in $\Omega$, there exists a positive measure
$\widehat{\sigma}$, supported in $\partial \Omega$, such that
$\|\sigma\|=\|\widehat{\sigma}\|$ and
\begin{equation}\label{balayage}
V^{\widehat{\sigma}}(z)-V^{\sigma}(z)\,=\,\const \,,\;
\text{qu.e.}\; z\notin \Omega \,,
\end{equation}
\noindent where $\dsty \const\,=\,0\,$ when $\Omega\,$ is a
bounded set, and a property is said to be satisfied for ``quasi-every" (qu.e.)
$z$ in a certain set, if it holds except for a possible subset of zero (logarithmic) capacity.
Then, $\dsty \widehat{\sigma}\,$ is said to be the
balayage of $\sigma\,$ from $\Omega\,$ onto $\partial \Omega\,.$

Now, we have the following:

\begin{lemma} \label{lem:measure}
The a priori complex measure \eqref{measures} is a unit positive
measure in $\Gamma_r$ \eqref{szegolevel}, for $0\leq r<+\infty$ . Moreover, $\mu_r$ is
the balayage of $\delta_0$ from $G_r$ onto $\Gamma_r$, where
$\delta_0$ denotes the Dirac Delta at $z=0$.
\end{lemma}

Now, for each $n\in \mathbb{N}\,,$ consider the ``pathological''
subset of negative integers $\dsty
\mathbb{S}_n\,=\,\{-n,-(n-1),\ldots,-2,-1\}\,.$ Hereafter, suppose that $\dsty \alpha_n \notin \mathbb{S}_n\,.$

Finally, denote by $\dsty \dist(\alpha_n, \mathbb{
S}_n)\,>\,0\,$ the minimal distance between $\alpha_n\,$ and the set $\mathbb{S}_n\,.$

\begin{theorem}\label{thm:main}
Consider a sequence of rescaled Laguerre polynomials
$\{L_n^{(\alpha_n)}(nz)\}_{n\in \mathbb{N}}$, such that
$\displaystyle \lim_{n\rightarrow
\infty}\frac{\alpha_n}{n}\,=\,-1\,$ and
\begin{equation}\label{limdist1}
\lim_{n\rightarrow\infty} [\dist(\alpha_n, \mathbb{S}_n)]^{1/n}\,=\,e^{-r}\,,
\end{equation}
\noindent for some $r\geq 0\,.$ Then, the contracted zeros of Laguerre polynomials asymptotically
follow the measure $d\mu_r$ in \eqref{measures} on the curve $\dsty \Gamma_r\,$ \eqref{szegolevel}. For $r=+\infty$, the limit measure is
$d\mu_{\infty}\,=\,\delta_0$.

\end{theorem}

\begin{remark}
The results above also hold when dealing with infinite
subsequences $\{L_n^{(\alpha_n)}(nz)\}_{n\in \Lambda}\,,\;\Lambda \subset \mathbb{N}\,.$
\end{remark}

\begin{figure}[htb] \label{fig:generic}
\centering \includegraphics[scale=0.7]{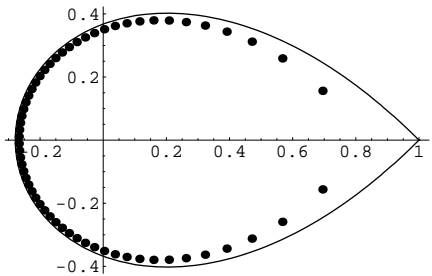} \caption{The
Szeg\H{o} curve and the zeros of $L_{60}^{(-60.1)}(60z)$.}
\end{figure}

\begin{remark}\label{rem:generic}
Observe that the case $\dsty r=0\,$ in Theorem \ref{thm:main} is
generic, because it takes place when parameters $\alpha_n$ do not
approach, or, at least, do not approach exponentially fast, the
set of integers $\dsty \mathbb{S}_n$ (see Figure \ref{fig:generic}). On the other hand, when
$r>0$ and, so, parameters approach the set of integers
$\dsty \mathbb{S}_n$ exponentially fast, the Szeg\H{o} curve
$\Gamma$ is replaced by a level curve $\Gamma_r$ which surrounds
$z=0$ and is strictly contained in the interior of $\Gamma$ (see
Figure 3). Finally, when $r=\infty$, i.e., when parameters
approach the set $\dsty \mathbb{S}_n$ faster than
exponentially, the limit measure reduces to a Dirac mass at $z=0$.
\end{remark}

\begin{figure}[htb] \label{fig:nongener}
\centering \includegraphics[scale=0.7]{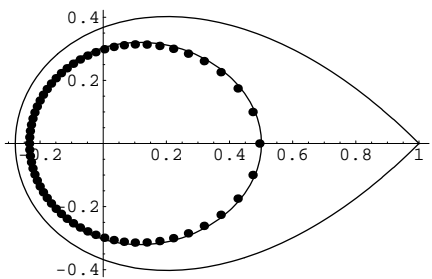} \caption{Zeros of
$L_{60}^{(-60+10^{-5})}(60z)$ and the curve $\Gamma_r$, for
$r=\frac{1}{12}\ln 10$.}
\end{figure}

\begin{remark}\label{rem:matching}
The weak asymptotics in the case $\dsty A=-1\,,$ characterized for
the set of measures \eqref{measures}
and the corresponding set of closed curves \eqref{szegolevel}, is the natural
matching between the solutions of the cases $\dsty A<-1\,$ (see
\cite{Kuijlaars01} and \cite{MR1858305}) and $\dsty -1<A<0\,$ (see
\cite{Kuijlaars04}). For those cases the following full set of
non-hermitian orthogonality relations for the Laguerre polynomials
with parameters $\dsty \alpha\in\mathbb{C}\,$ was used:
\begin{equation*}\label{nonhermit}
\int_{\Sigma}\,L_n^{(\alpha)}(z)\,z^k\,z^{\alpha}\,e^{-z}\,dz\,=\,0\,,\;
k=0,\ldots,n-1\,,
\end{equation*}
\noindent where $\Sigma\,$ is any unbounded contour in
$\mathbb{C}\setminus [0,\infty)\,,$ connecting $+\infty +iy\,$ and
$+\infty -iy\,,$ for some $y>0\,,$ and the branch in
$z^{\alpha}\,$ is taken with the cut along the positive real axis
(see \cite[Lemma 2.1]{Kuijlaars01}). In \cite{MR1858305} this full
set of orthogonality relations allowed to apply seminal results by
H. Stahl \cite{Stahl:86} and A. Gonchar and E. A. Rakhmanov
\cite{Gonchar:89} on the asymptotic behavior of complex orthogonal
polynomials. Indeed, it was proved that zeros of the rescaled
Laguerre polynomials accumulate on a closed contour $C\,$ in
$\mathbb{C}\setminus [0,\infty)\,$ which is ``symmetric'' (in the
``Stahl's sense'', see \cite{Stahl:85}-\cite{Stahl:86}) with respect
to the external field $\varphi(z)\,=\,\frac{1}{2}\,(-A\log|z|+\Re
z)\,,$ and that they asymptotically follow the equilibrium
distribution on $C\,$ in presence of the external field
$\varphi\,$. In the proof of the main result in this paper,
it will be showed that the zeros of the rescaled Laguerre polynomials in the present case
also asymptotically follow the equilibrium distribution of $\Gamma_r\,$ in presence of
the external field $\varphi\,$ above (for $A=-1\,$), $\Gamma_r\,$
being a symmetric contour with respect to this external field.
That is, although the theorems by H. Stahl and A. Gonchar-E. A.
Rakhmanov cannot be applied in this case since the complement to the support is
disconnected, the conclusions still hold.
\end{remark}

\section{On asymptotically extremal polynomials}

Throughout this section, some topics in potential theory which are
needed for the proof of our main result will be recalled. For more
details the reader can consult the monography \cite{Saff:97}.

First, let us precise the notion of admissible weights.

\begin{definition} \label{def:admissible}
Given a closed set $\Sigma \subset \mathbb{C}$, we say that a
function $\omega : \Sigma \longrightarrow [0,\infty)$ is an
admissible weight on $\Sigma$ if the following conditions are
satisfied (see \cite[Def.I.1.1]{Saff:97}):

\begin{itemize}

\item [(a)] $\omega$ is upper semi-continuous;

\item [(b)] the set $\{z\in \Sigma : \omega (z)>0\}$ has positive
(logarithmic) capacity;

\item [(c)] if $\Sigma$ is unbounded, then $\displaystyle
\lim_{|z|\rightarrow \infty ,\,z\in \Sigma} |z|\omega (z)=0$.

\end{itemize}

\end{definition}

Given such an admissible weight $\omega$ in the closed set
$\Sigma$, and setting $\varphi (z)= -\log \omega (z)$, we know
(see e.g. \cite[Ch.I]{Saff:97}) that there exists a unique measure
$\mu_{\omega}$, with (compact) support in $\Sigma$, for which the
infimum of the weighted (logarithmic) energy
\begin{equation*}\label{weightedenergy}
I_{\omega}(\mu)\,= \,-\int\int \log |z-x|d\mu (z)d\mu
(x)\,+\,2\int \varphi (x)d\mu (x)
\end{equation*}
\noindent is attained. Moreover, setting $\displaystyle F_{\omega}
= I_{\omega}(\mu_{\omega})-\int \varphi d\mu_{\omega}\,,$ which is called the modified Robin constant, we have
the following property, which uniquely characterizes the extremal
measure $\mu_{\omega}$:
\begin{equation*}\label{equilibrium1}
V^{\mu_{\omega}}(z)+\varphi (z)\,\begin{cases} =\,F_{\omega}\,,\,
& qu.e.\; z\in \supp \mu_{\omega}\,, \\
\geq \,F_{\omega}\,,\, & qu.e.\; z\in \Sigma \,,
\end{cases}
\end{equation*}
\noindent where for a measure $\sigma$, $V^{\sigma}$ denotes its
logarithmic potential, that is,
$$
V^{\sigma}(z)\,= \,-\int \log |z-x|\,d\sigma (x)\,.
$$
Now, let $\Sigma$ be a closed set and $\omega$ an admissible
weight on $\Sigma$. Then, a sequence of monic polynomials
$\{p_n\}_{n\in \mathbb{N}}$ is said to be asymptotically extremal
with respect to the weight $\omega$ if it holds (see
\cite{Saff:97}):
\begin{equation}\label{asympextr1}
\lim_{n\rightarrow
\infty}\|\omega^np_n\|^{1/n}_{\Sigma}\,=\,\exp
(-F_{\omega})\,,
\end{equation}
\noindent where, as usual, $\dsty \|\cdot\|_K\,$ denotes the sup-norm in the set $K\,.$
The study of weighted polynomials of the form $\displaystyle
\omega (z)^nP_n(z)$ has applications to many problems in
approximation theory (see e.g. the monographies \cite{Saff:97} and
\cite{To94}). It is well known that if for each $n\in \mathbb{N}$,
$T_n^{\omega}$ is the $n$-th (weighted) Chebyshev polynomial with
respect to the weight $\omega ^n$, that is, if it is the (unique) monic
polynomial of degree $n$ for which the infimum
$$
t_n^{\omega}\,=\,\inf \{\|\omega
^nP\|_{\Sigma}\,,\,P(z)=z^n+\ldots\}
$$
\noindent is attained, then the sequence $\{T_n^{\omega}\}$
satisfies the asymptotic behavior given in \eqref{asympextr1} (see
\cite[Ch.III]{Saff:97}).

Under mild conditions on the weight $\omega$, in
\cite[Ch.III]{Saff:97} it is shown that the zeros of such
sequences of polynomials asymptotically follow the equilibrium
measure $\mu_{\omega}$, in the sense of the weak-* convergence.
Indeed, we have the following result (see
\cite[Th.III.4.1]{Saff:97} or the previous paper \cite{MhaSaf91}):

\begin{theorem} \label{thm:SaTo}

Let $\omega$ be an admissible weight such that the support of the
corresponding equilibrium measure $\mu_{\omega}$, $S_{\omega}$,
has zero Lebesgue planar measure. Let $\{p_n\}_{n\in \mathbb{N}}$ be a sequence of
monic polynomials of respective degrees $n=1,2,\ldots$ satisfying:

\begin{equation}\label{extremal1}
\lim_{n\rightarrow \infty}\left\|\omega^n
p_n\right\|_{S_{\omega}}^{1/n}\,=\,exp\,(-F_{\omega})\,,
\end{equation}

\noindent where $F_{\omega}$ denotes the modified (by the external
field $\varphi = -\ln \omega$) Robin constant. Then, the following
statements are equivalent:

\begin{itemize}

\item [(a)] $\nu (p_n) \longrightarrow \mu_{\omega}$ in the weak-* sense, where $\nu
(p_n)$ denotes the unit zero counting measure associated to $p_n$,
that is, $\dsty d\nu
(p_n)=\frac{1}{n}\,\sum_{p_n(\zeta)=0}\delta_{\zeta}\,.$

\item [(b)] For each bounded component $\textit{R}$ of
$\mathbb{C}\setminus S_{\omega}$ and each infinite sequence
$\textit{N} \subset \mathbb{N}$, there exist $z_0 \in \textit{R}$
and $\textit{N}_1 \subset \textit{N}$ such that
\begin{equation}\label{extremal2}
\lim_{n\rightarrow \infty,n\in \textit{N}_1}\left|
p_n(z_0)\right|^{1/n}\,=\,exp\,(-V^{\mu_{\omega}}(z_0))\,.
\end{equation}

\end{itemize}

\end{theorem}

\begin{remark}\label{rem:GaKu}

In \cite[Theorem 5]{GaKu97}, condition \eqref{extremal1} is replaced by the weaker condition:
\begin{equation}\label{extremalGaKu}
\limsup_{n\rightarrow \infty}\,\omega (z) |p_n(z)|^{1/n}\,\leq
\,exp\,(-F_{\omega})\,,\; \text{qu.e.} \; z \in S_{\omega}.
\end{equation}
\end{remark}

\begin{remark}\label{rem:externalfield}
It is clear that the balayage of a measure (see \eqref{balayage}) is a very particular case of equilibrium measure in an external field. Since Lemma \ref{lem:measure}  says that measure $\dsty \mu_r\,$ is the balayage of $\dsty
\delta_0\,$ from $\dsty G_r\,$ onto its boundary $\dsty \Gamma_r\,,$ it means that
\begin{equation}\label{balay}
V^{\mu_r}(z)\,=\,-\log|z|\,,\;z\in \Gamma_r\,.
\end{equation}
Taking into account the expression of $\dsty \Gamma_r\,,$ \eqref{balay} implies both
\begin{equation}\label{externalRe}
V^{\mu_r}(z)\,+\,\Re z\,=\,r+1\,,\;z\in \Gamma_r\,,
\end{equation}
\noindent and
\begin{equation}\label{externalStahl}
V^{\mu_r}(z)\,+\,\varphi(z)\,=\,\frac{r+1}{2}\,,\;z\in \Gamma_r\,,
\end{equation}
\noindent where
\begin{equation}\label{varphi}
\varphi(z)\,=\,\frac{1}{2}\,\left(\log|z|+\Re z\right)\,,
\end{equation}
\noindent (see Remark \ref{rem:matching} above).
\end{remark}

For the proof of Theorem \ref{thm:main}, taking into account Theorem \ref{thm:SaTo}, it
will be proved that the rescaled Laguerre polynomials are
asymptotically extremal with respect to the weight $\dsty \omega=e^{-\varphi}\,$ in the compact set given by the closed
contour $\Gamma_r\,$ (using \eqref{extremalGaKu}), along with the fact that they satisfy the local behavior \eqref{extremal2}.

\section{Proofs}

\subsection{Proof of Lemma \ref{lem:measure}}

Take into account that the level curves $\Gamma_r$, for $0\leq
r<\infty$, given by \eqref{szegolevel} are, in fact, trajectories
of the quadratic differential (see e.g. \cite{Strebel})
$$
-\,\frac{(z-1)^2}{z^2}\,\,(dz)^2\,\,,
$$
or, what is the same, $\Gamma_r$ may be defined in the form:
\begin{equation}\label{traject}
\Gamma_r\,=\,\left\{z\in
\mathbb{C}/\,\Re\,\int_1^z\left(1-\frac{1}{t}\right)\,dt\,=\,r\right\}\,.
\end{equation}
Expression \eqref{traject} shows that \eqref{measures} is
real-valued in $\Gamma_r$ and does not change its sign. Moreover,
by a straightforward application of the Cauchy theorem, we have
that
$$
\mu_r(\Gamma_r)\,=\,\int_{\Gamma_r}d\mu_r(t)\,=\,1\,.
$$
Now, we will prove that $\mu_r$ is the balayage of $\delta_0$ from
$G_r=\Int(\Gamma_r)$ onto $\Gamma_r$.

To this end, consider the function $\displaystyle
\phi(z)=ze^{1-z}$. It is easy to see that $\phi$ conformally maps
$G_r$ onto the disk $\mathbb{D}_r=\{w\in
\mathbb{C}/|w|<r\}\,,\;0\leq r<\infty\,,$ in the $w$-plane (see
\cite{szego:1924} and \cite{PrVa97}). Thus, from \eqref{measures},
we have:
$$
d\mu_r(z)\,=\,\frac{1}{2\pi
i}\,\left(\frac{1}{z}-1\right)dz\,=\,\frac{1}{2\pi
i}\,\frac{\phi'(z)}{\phi(z)}\,dz\,=\,\frac{1}{2\pi
i}\,\frac{dw}{w}\,=\,\frac{d\theta}{2\pi}\,,
$$
where $\dsty w=r\,e^{i\theta}=\phi (z)\,,$ and $\dsty
z\in\Gamma_r\,.$ Therefore, \eqref{measures} is the preimage of
the normalized arc-length measure on the circle $\mathbb{T}_r\,=\,\partial \mathbb{D}_r\,$
under the mapping $w=\phi(z)$, that is, the harmonic measure at
$z=0$ with respect to the domain $G_r$. But this fact implies that
\eqref{measures} is the balayage of $\delta_0$ from $G_r$ onto
$\Gamma_r$ (see \cite[p. 222]{Landkof}).

\subsection{Proof of Theorem \ref{thm:main}}

In Remark \ref{rem:externalfield}, it was shown that $\dsty \mu_r\,$ is also the equilibrium measure in $\dsty \Gamma_r\,$ in the external field
$\dsty \varphi\,$ \eqref{varphi}.

Moreover, \eqref{externalStahl} shows that the
corresponding modified Robin constant is given by:
\begin{equation}\label{Robin}
F_{\omega}\,=\,\frac{r+1}{2}\,.
\end{equation}
On the other hand, the function $\dsty g(z)\,=\,V^{\mu_r}(z)+\Re z\,$ is harmonic in $\overline{G_r}\,$ and, by \eqref{externalRe}, $\dsty g(z)\,\equiv\,r+1\,,\;z\in \Gamma_r\,.$ Then,
it yields that $\dsty g(z)\,\equiv\,r+1\,,\;z\in \overline{G_r}\,.$ In particular,
\begin{equation}\label{zero}
V^{\mu_r}(0)\,=\,r+1\,.
\end{equation}
From \eqref{Robin}, in order to prove \eqref{extremalGaKu} we need to
show that
\begin{equation*}
\limsup_{n\rightarrow \infty}\;\omega (z) |p_n(z)|^{1/n}\,\leq
\,e^{-\frac{r+1}{2}}\,,\; \text{qu.e.}\; z\in \Gamma_r\,,
\end{equation*}
\noindent for the monic polynomial $\displaystyle
p_n(z)\,=\,\widehat{L}_n^{(\alpha_n)}(nz)$ and the weight
$\displaystyle \omega (z)= e^{-\varphi
(z)}\,,$ which taking into
account the expression of $\Gamma_r$, is equivalent to prove:
\begin{equation} \label{extremal11}
\limsup_{n\rightarrow \infty}\; e^{-\Re z}
|p_n(z)|^{1/n}\,\leq \,e^{-(r+1)}\,,\;\text{qu.e.}\; z\in \Gamma_r.
\end{equation}
Now, since by \eqref{explLag} $\displaystyle
L_n^{(\alpha_n)}(nz)\,=\,l_n^{\alpha_n}z^n+\ldots$, with
\begin{equation}\label{director}
l_n^{\alpha_n}=(-1)^n\frac{n^n}{n!}\,,
\end{equation}
we have that \eqref{extremal11} is equivalent to
\begin{equation} \label{extremal111}
\limsup_{n\rightarrow \infty}\; e^{-\Re z}
|L_n^{(\alpha_n)}(nz)|^{1/n}\,\leq \,e^{-r}\,,\;\text{qu.e.}\; z\in \Gamma_r.
\end{equation}
In addition, we should prove that there exists a point $z_0\in
G_r$ for which \eqref{extremal2} is attained. Thus, choosing
$z_0=0$, and taking into account \eqref{zero}, it is enough to show that
\begin{equation*}
\lim_{n\rightarrow \infty}\left|
p_n(0)\right|^{1/n}\,=\,e^{-(r+1)}\,,
\end{equation*}
\noindent or what is the same, by \eqref{director},
\begin{equation} \label{extremal22}
\lim_{n\rightarrow
\infty}\left|L_n^{(\alpha_n)}(0)\right|^{1/n}\,=\,e^{-r}\,.
\end{equation}
Now, we are going to prove \eqref{extremal22} and
\eqref{extremal111} under the conditions in Theorem \ref{thm:main}.

\subsubsection{Proof of \eqref{extremal22}}

Take into account that by
\eqref{explLag},
\begin{equation*}
L_n^{(\alpha_n)}(0)\,=\,\binom{n+\alpha_n}{n}\,,
\end{equation*}
\noindent and let $\dsty h_n\in \{1,2,\ldots n\}$ be such that
    $$
\dist(\alpha_n, \mathbb{S}_n)\,=\,|\alpha_n+h_n|\,.
    $$
Thus, by \eqref{limdist1}, we have:
    $$
\lim_{n\rightarrow \infty}|\alpha_n+h_n|^{1/n}\,=\,e^{-r}\,,
    $$
\noindent and, therefore, to prove \eqref{extremal22} it should be
satisfied:
\begin{equation}\label{identityf}
\dsty \lim_{n\rightarrow\infty} \left(\frac{\left|(n + \alpha_n)(n
+ \alpha_n - 1)\ldots(1 + \alpha_n)\right|} {n! |\alpha_n +
h_n|}\right)^{\frac{1}{n}} = 1\, .
\end{equation}

Let us suppose, first, that $\alpha_n \geq  -n - \frac{1}{2}\,.$ Then, it is easy to see that
$$
\left|(n + \alpha_n)(n + \alpha_n - 1)\ldots(1 + \alpha_n)\right|
\,=\,\left|\alpha_n+h_n\right|\,\prod_{k=1}^{n-h_n}\,\left|\alpha_n+h_n+k\right|\,\prod_{k=1}^{h_n-1}\,
\left|\alpha_n+h_n-k\right|\,,
$$
\noindent and taking into account that $\dsty \frac{2k-1}{2} \leq
\left|\alpha_n + h_n\, \pm \,k\right|\leq \frac{2k+1}{2}\,,$ for
any integer $k\geq 1\,,$ it yields
$$
\begin{array}{rcccl}
\dsty \prod_{k=1}^{h_n-1} \frac{2k - 1}{2}\prod_{k=1}^{n-h_n}
\frac{2k - 1}{2}& \leq &\dsty \frac{\left|(n + \alpha_n)(n +
\alpha_n - 1)\ldots(1 + \alpha_n)\right|}{|\alpha_n + h_n|} & \leq
& \dsty \prod_{k=1}^{h_n-1} \frac{2k + 1}{2}\prod_{k=1}^{n-h_n}
\frac{2k + 1}{2}\,.
\end{array}
$$
Thus, denoting $a_l = \dsty \prod_{k=1}^{l} \frac{2k + 1}{2} =
\frac{2(l+1)!}{2^{2(l+1)-1}(l+1)!}\,$, $\dsty l\geq 1\,,$ and $\dsty a_0=1\,, a_{-1}=2\,.$ Then,
$$
\begin{array}{rcccl}
\dsty \frac{1}{2^2} a_{h_n-2} a_{n-h_n-1} & \leq & \dsty
\frac{\left|(n + \alpha_n)(n + \alpha_n - 1)\ldots(1 +
\alpha_n)\right|}{|\alpha_n + h_n|} & \leq  & \dsty a_{h_n-1}
a_{n-h_n}\,,\, 1  \leq h_n \leq n\,.
\end{array}
$$
On the other hand, if  $\alpha_n < -n - \frac{1}{2}\,$ (and thus, $h_n=n\,$),
$$
\begin{array}{rcccccl}
\dsty
\left|(n + \alpha_n - 1)(n + \alpha_n - 2)\ldots(1 +
\alpha_n)\right| & \leq  & \dsty \prod_{k=1}^{n-1}(-\alpha_n  - k ) & = &  \dsty \frac{\Gamma(-\alpha_n)}{\Gamma(-(\alpha_n + n -1))}& = & \dsty \frac{\Gamma(-\alpha_n) }{\Gamma(\delta_n + 1)}\,,
\end{array}
$$
\noindent from which it yields
$$
\begin{array}{rcccl}
\dsty  a_{n-1} \;\; \leq \;\; \dsty
\left|(n + \alpha_n - 1)\ldots(1 +
\alpha_n)\right| \;\; \leq  \;\;\frac{\Gamma(-\alpha_n))}{\Gamma(\delta_n + 1)}\,.
\end{array}
$$
Now, since
$$\lim_{n\rightarrow\infty} \left(\frac{a_{h_n-1 - s} a_{n-h_n-s}}{n!}\right)^{\frac{1}{n}} = 1\,,\, s = 0,1 \,,$$
\noindent and
$$\lim_{n\rightarrow \infty}\left(\dsty\frac{\Gamma(-\alpha_n)}{n!\,\Gamma(\delta_n + 1)}\right)^{\frac{1}{n}} = 1\,,$$
\noindent then \eqref{identityf} follows.

\subsubsection{Proof of \eqref{extremal111} }

Let us denote
\begin{equation}\label{ceiling}
k_n\,=\,\min\,\left([-\alpha_n],n\right)\,,\;\;\alpha_n\,=\,-k_n-\delta_n\,,\;\delta_n>0\,,
\end{equation}
\noindent where, as usual, $\dsty [\cdot]\,$ denotes the integer
part of a real number. It is clear that $\dsty -k_n\in
\mathbb{S}_n\,$ and if $\dsty k_n<n\,,$ then $\dsty
0<\delta_n<1.$

It also holds
    $$
\dist(\alpha_n, \mathbb{S}_n)\,=\,\begin{cases}\,\delta_n\,,\; & \text{if}\;\alpha_n<-n\,, \\ \min (\delta_n,1-\delta_n)\,,\; & \text{if}\;\alpha_n>-n\,.
\end{cases}
    $$
In order to prove \eqref{extremal111}, the
following integral representation will be used (see \cite[formula
(6.2.22)]{Askey}):
\begin{equation}\label{repAskeyun}
e^{-x}L_n^{(\alpha)}(x)\,=\,\frac{1}{\Gamma(\beta-\alpha)}\,\int_x^{\infty}\,(t-x)^{\beta-\alpha-1}e^{-t}L_n^{(\beta)}(t)\,dt\,,
\end{equation}
\noindent where $\beta>\alpha\,$ and the path of integration is
any simple smooth path connecting $x\in \mathbb{C}\,$ with
$+\infty\,.$ Thus, setting $\dsty \beta=-k_n\,$ and $\,\alpha=\alpha_n\,$ in
\eqref{repAskeyun} and taking into account \eqref{ceiling}, we
have:
\begin{equation*}\label{repAskey1un}
e^{-x}\,L_n^{(\alpha_n)}(x)\,=\,\frac{1}{\Gamma(\delta_n)}\,\int_x^{\infty}\,(t-x)^{\delta_n-1}\,
e^{-t}\,L_n^{(-k_n)}(t)\,dt\,,
\end{equation*}
\noindent or what is the same, after some calculations,
\begin{equation}\label{repAskey2un}
e^{-nx}\,L_n^{(\alpha_n)}(nx)\,=\,\frac{n^{\delta_n}}{\Gamma(\delta_n)}\,\int_x^{\infty}\,(t-x)^{\delta_n-1}\,
e^{-nt}\,L_n^{(-k_n)}(nt)\,dt\,.
\end{equation}
Now, since $\dsty k_n\in\{1,\ldots,n\}\,,$ making use of
\eqref{expneg}, \eqref{repAskey2un} may be written in the form:
\begin{equation}\label{repAskey3un}
e^{-nx}\,L_n^{(\alpha_n)}(nx)\,=\,(-1)^{k_n}\,\frac{n^{\delta_n+k_n}\,(n-k_n)!}{n!\,
\Gamma(\delta_n)}\,\int_x^{\infty}\,(t-x)^{\delta_n-1}\,t^{k_n}
e^{-nt}\,L_{n-k_n}^{(k_n)}(nt)\,dt\,.
\end{equation}

On the other hand, taking into account the Rodrigues formula
\eqref{RodrLag}, \eqref{repAskey3un} yields:
\begin{equation*}\label{repAskey4un}
\begin{split}
e^{-nx}\,L_n^{(\alpha_n)}(nx)\,=\, & (-1)^{n}\,e^{-n}\,\frac{n^{\delta_n+k_n}}{n!\,
\Gamma(\delta_n)}\,\int_x^{\infty}\,(t-x)^{\delta_n-1}\,\left[\phi(t)^n\right]^{(n-k_n)}\,dt\, \\
= \,&   \Delta_n\,F_n(x)\,,
\end{split}
\end{equation*}
\noindent where, as above, $\dsty \phi(t)=te^{1-t}\,$ and $$\dsty F_n(x)\,=\,\int_x^{\infty}\,(t-x)^{\delta_n-1}\,\left[\phi(t)^n\right]^{(n-k_n)}\,dt\,.$$
Let us denote by $x_0\,=\,x_0(r)\,$ the unique point where the curve $\Gamma_r$ meets
the positive real semiaxis.  Now, using the freedom in the choice of the path of integration, it will consists of two arcs: the first goes from $x$ to $x_0$ through
the curve $\Gamma_r\,$ (by the shortest way), and the corresponding integral will be
denoted by $\dsty G_n(x)$; the second goes from $x_0$ to
$\infty\,$ along the positive real semiaxis, and we will denote this integral
by $\dsty H_n(x)\,.$ Thus, $\dsty F_n(x)=G_n(x)+H_n(x)\,.$


We are going to estimate $G_n(x)\,,$ for $x\in \Gamma_r\setminus \{x_0\}\,.$

Suppose first that $k_n=n\,,$ and hence,
$$\dsty G_n(x)\,=\,\int_x^{x_0}\,(t-x)^{\delta_n-1}\,\phi(t)^n\,dt\,.$$
For it, consider the natural arc-length
parametrization: $t=t(s)\,,$ so that $t(0)=x\,$ and
$t(s_0)=x_0\,,$ for some positive real number $s_0\,.$ In addition, recall that $\dsty
|\phi(t)|= e^{-r}\,,$ for $\dsty t\in \Gamma_r\,.$ Since the path
of integration  is a smooth recitifiable Jordan arc (even for the
case when $r=0\,,$ since the path is entirely contained in the
upper, or lower, half of $\Gamma_0 = \Gamma\,$), we have
\begin{equation}\label{path1}
\begin{split}
|G_n(x)|= & \left|\int_0^{s_0}\,(t(s)-t(0))^{\delta_n-1}\,(\phi(t(s)))^n\,t'(s)\,ds\right| \\
& \leq \,\|\phi\|_{\Gamma_r}^n\,
\int_0^{s_0}\,|t(s)-t(0)|^{\delta_n-1}|t'(s)|ds\, \\
& \leq \,e^{-rn}
\int_0^{s_0}\,|t(s)-t(0)|^{\delta_n-1}|t'(s)|ds\,.
\end{split}
\end{equation}
Now, take into account that there exist two positive constants $k,C\,,$ such that $\dsty k\leq
|t'(s)|\leq C\,,\; s\in [0,s_0]\,,$ and set $A_n\,=\begin{cases} C^{\delta_n}\,,\; & \text{if}\; \delta_n\geq 1 \,, \\ C\,k^{\delta_n-1}\,,\; & \text{if}\; 0<\delta_n<1 \,. \end{cases}$ Then, by classical mean value theorem, \eqref{path1} implies:
\begin{equation}\label{boundGn}
|G_n(x)|\,\leq \,A_n\,e^{-rn}\,\frac{s_0^{\delta_n}}{\delta_n}\,,
\end{equation}
\noindent where $\dsty \lim_{n\rightarrow \infty}\,A_n^{1/n}\,=\,1\,.$
On the other hand, when $k_n<n\,,$ it follows
$$\dsty G_n(x)\,=\,\int_x^{x_0}\,(t-x)^{\delta_n-1}\,\left[(\phi(t))^n\right]^{(n-k_n)}\,dt\,.$$
Proceeding analogously as above, it holds
$$ |G_n(x)|\,\leq \,A_n\,\left\|[\phi^n]^{(n-k_n)}\right\|_{\Gamma_r}\,\frac{s_0^{\delta_n}}{\delta_n}\,,$$
and thus, by applying the Cauchy integral formula in an arbitrarilly small
circle around $t\,,$ we have for $t$ in the segment of curve
$\dsty \Gamma_r\,$ connecting $x$ to $x_0\,,$
\begin{equation*}\label{phideriv}
\begin{split}
\left|\left[\phi(t)^n\right]^{(n-k_n)}\right|\,\leq & \,(n-k_n)!\,\epsilon^{-n+k_n}\,e^{\epsilon
n}\,(|\phi(t)|+\epsilon\,e^2)^n\,= \\
&  \,(n-k_n)!\,\epsilon^{-n+k_n}\,e^{\epsilon
n}\,(e^{-r}+\epsilon\,e^2)^n\,,
\end{split}
\end{equation*}
\noindent for $\dsty \epsilon>0\,$ arbitrarily small. Hence,
\begin{equation}\label{Gn2}
|G_n(x)|\,\leq \,A_n\,(n-k_n)!\,\epsilon^{-n+k_n}\,e^{n\varepsilon}\,\left(e^{-r}+\varepsilon\,e\right)^n\,\frac{s_0^{\delta_n}}{\delta_n}\,,
\end{equation}
\noindent for $\varepsilon >0\,.$
Since $\dsty \lim_{n\rightarrow \infty}\,\frac{k_n}{n}\,=\,1\,,$ we have,
\begin{equation}\label{rootGn2}
\lim_{n\rightarrow \infty}\left((n-k_n)!\,\epsilon^{-n+k_n}\,e^{n\varepsilon}\,\left(e^{-r}+\varepsilon\,e\right)^n\,\right)^{1/n}\,=
\,e^{\varepsilon}\,\left(e^{-r}+\varepsilon\,e\right)\,,
\end{equation}
\noindent for $\varepsilon >0\,.$
Observe that taking $\dsty \varepsilon \rightarrow 0^+\,,$ from \eqref{rootGn2}, \eqref{Gn2} agrees with \eqref{boundGn}.
Taking into account that $\dsty \lim_{n\rightarrow \infty}\,\frac{\alpha_n}{n}\,=\,-1\,,$ we have
\begin{equation}\label{identity}
\lim_{n\rightarrow\infty}\,\left(|\Delta_n|\,\frac{s_0^{\delta_n}}{\delta_n}\right)^{1/n}\,=\,\lim_{n\rightarrow
\infty}\,\frac{n^{\delta_n/n}}{\Gamma(1+\delta_n)^{1/n}}\,=\,1\,,
\end{equation}
\noindent where Stirling formula has been used when $\delta_n\,$ is unbounded (recall that $\delta_n=o(n)\,$). Therefore, by \eqref{boundGn}-\eqref{identity}, it yields
\begin{equation}\label{DeltanGn}
\limsup_{n\rightarrow\infty}\,\left(|\Delta_n\,G_n(x)|\right)^{1/n}\,\leq \,e^{-r}\,,\;x\in\Gamma_r\setminus\{x_0\}\,,
\end{equation}
\noindent after taking limits when $\dsty \varepsilon \rightarrow 0^+\,,$ if necessary. Note that in this part of the proof \eqref{limdist1} has not been used.

Now, we are concerned with $H_n(x)\,.$ As above, suppose first that $k_n=n\,,$ and thus,
\begin{equation*}
H_n(x)\,=\,\int_{x_0}^{\infty}\,(t-x)^{\delta_n-1}\,\phi(t)^n\,dt\,,
\end{equation*}
\noindent where now the path of integration is contained in the positive real semiaxis. Then, we have
\begin{equation}\label{Hn88}
H_n(x)\,=\,e^n\,\int_{x_0}^{\infty}\,\left(1-\frac{x}{t}\right)^{\delta_n-1}\,t^{n+\delta_n-1}\,e^{-(n-1)t}\,e^{-t}\,dt\,.
\end{equation}
Taking into account that there exist two positive constants $M,N\,,$ such that $\dsty M\leq
\left|1-\frac{x}{t}\right|\leq N\,,\; t\in [x_0,\infty)\,,$ and setting $B_n\,=\begin{cases} N^{\delta_n-1}\,,\; & \text{if}\; \delta_n\geq 1\,, \\ M^{\delta_n-1}\,,\; & \text{if}\; 0<\delta_n<1\,, \end{cases}$ then,
\begin{equation*}
|H_n(x)|\, \leq \,e^n\,B_n\,\|h\|_{[0,+\infty)}\,\int_{x_0}^{\infty}\,e^{-t}\,dt\,
 \leq \,e^n\,B_n\,\|h\|_{[0,+\infty)}\,,
\end{equation*}
\noindent where $\dsty \lim_{n\rightarrow \infty}\,B_n^{1/n}\,=\,1\,,$ $\dsty h(t)\,=\,t^{n+\delta_n-1}\,e^{-(n-1)t}\,$ and it is not hard to see that
$$\dsty \|h\|_{[0,+\infty)}\,= h\left(\frac{n+\delta_n-1}{n-1}\right)\,=\,\left(\frac{n+\delta_n-1}{n-1}\right)^{n-\delta_n-1}\,e^{-(n+\delta_n-1)}\,.$$ Hence,
\begin{equation*}\label{Hn11}
|H_n(x)|\,  \leq \,B_n\,\left(\frac{n+\delta_n-1}{n-1}\right)^{n-\delta_n-1}\,e^{-(\delta_n-1)}\,.
\end{equation*}
Therefore,
\begin{equation}\label{Hn22}
\begin{split}
|\Delta_n\,H_n(x)|\,  & \leq \,\frac{e^{-n}n^{n+\delta_n}}{n!\Gamma(1+\delta_n)}\,B_n\,\left(\frac{n+\delta_n-1}{n-1}\right)^{n-\delta_n-1}\,e^{-(\delta_n-1)}\,\delta_n\, \\[.5cm]
& \leq C_n\,\delta_n\,=\,C_n\,\dist(\alpha_n,\mathbb{S}_n)\,,
\end{split}
\end{equation}
\noindent where $\dsty \lim_{n\rightarrow\infty}\,C_n^{1/n}\,=\,1\,.$

On the other hand, when $k_n<n\,,$ we have
$$
H_n(x)\,=\,\int_{x_0}^{\infty}\,(t-x)^{\delta_n-1}\,[\phi(t)^n]^{(n-k_n)}\,dt\,,
$$
\noindent and integrating by parts, it yields
\begin{equation}\label{parts}
H_n(x)\,=\, (x_0-x)^{\delta_n-1}\,[\phi(t)^n]^{(n-k_n-1)}_{t=x_0}\,+\,  (1-\delta_n)\int_{x_0}^{\infty}\,(t-x)^{\delta_n-2}\,[\phi(t)^n]^{(n-k_n-1)}\,dt\,.
\end{equation}
Now, applying again the Cauchy integral formula for $t\in
[x_0,\infty)\subset\mathbb{R}^+\,,$ it holds:
\begin{equation}\label{phiderivh}
\left|\left[(\phi(t))^n\right]^{(l)}\right|\,\leq\,l!\,\epsilon^{-l}\,e^{2\epsilon
n}\,\phi(t+\varepsilon)^n\,,
\end{equation}
\noindent for arbitrarily small $\epsilon>0\,.$

Then, taking into account \eqref{parts}-\eqref{phiderivh} and setting $$ D_n\,=\,(n-k_n -1)!\, \epsilon^{-n+k_n+1}\,e^{2\epsilon
n}\,|x_0-x|^{\delta_n-1}\,,$$ we have
\begin{equation*}
\begin{split}
|H_n(x)|\,  \leq \, D_n\,\left(\phi(x_0+\varepsilon)^n \,+\, (1-\delta_n)\,\int_{x_0}^{\infty}\,|t-x|^{-1}\,\phi(t+\varepsilon)^n\,dt\right) \\
 \leq \, D_n\,\left(\phi(x_0+\varepsilon)^n \,+\, (1-\delta_n)\,\int_{x_0}^{\infty}\,\left|1-\frac{x-\varepsilon}{t}\right|^{-1}\,t^{-1}\,\phi(t)^n\,dt\right)\,.
\end{split}
\end{equation*}
Finally, we can bound the integral above as in (\ref{Hn88}), which yields
\begin{equation*}\label{Hn12}
|H_n(x)|\,  \leq \, D_n\,\left(\phi(x_0+\varepsilon)^n\right. \,+
 \left. (1-\delta_n)\,\widetilde{M}^{-1}\,e^{-1}\right)\,,
\end{equation*}
\noindent where we denote by $\dsty \widetilde{M}\,$ the lower bound of the function $\dsty |1-\frac{x-\varepsilon}{t}|\,,\; t\in [x_0,\infty)\,.$

Therefore,
\begin{equation}\label{Hn13}
\begin{split}
|\Delta_n\,H_n(x)|\,& \leq \,\frac{e^{-n}n^{k_n+\delta_n}}{n!\Gamma (1+\delta_n)}\, D_n\,\left(\phi(x_0+\varepsilon)^n \,\delta_n\, +\, \delta_n(1-\delta_n)\,\widetilde{M}^{-1}\,e^{-1}\right) \\[.25cm]
& \leq R_n\,\delta_n\,\phi(x_0+\varepsilon)^n\,+\,S_n\,\delta_n (1-\delta_n)\,,
\end{split}
\end{equation}
\noindent where $\dsty \lim_{n\rightarrow\infty}\,R_n^{1/n}\,=\,\lim_{n\rightarrow\infty}\,S_n^{1/n}\,=\,1\,.$ Taking into account \eqref{limdist1}, we have
\begin{equation}\label{equivalence}
\lim_{n\rightarrow\infty}\,[\dist(\alpha_n,\mathbb{S}_n)]^{1/n}\,=\,\lim_{n\rightarrow\infty}\,[\delta_n (1-\delta_n)]^{1/n}\,=\,e^{-r}\,.
\end{equation}
Now, from \eqref{Hn22}, \eqref{Hn13} and \eqref{equivalence}, it follows
\begin{equation}\label{DeltanHn}
\limsup_{n\rightarrow\infty}\,\left(|\Delta_n\,H_n(x)|\right)^{1/n}\,\leq \,e^{-r}\,,\;x\in\Gamma_r\setminus\{x_0\}\,,
\end{equation}
\noindent after taking limits when $\dsty \varepsilon \rightarrow 0^+\,,$ if necessary. Thus, from \eqref{DeltanGn} and \eqref{DeltanHn}, it yields
\begin{equation*}\label{DeltanFn}
\limsup_{n\rightarrow\infty}\,\left(|\Delta_n\,F_n(x)|\right)^{1/n}\,\leq \,e^{-r}\,,\;x\in\Gamma_r\setminus\{x_0\}\,.
\end{equation*}

It only remains to consider the limit case $\dsty r=\infty \,,$ which occurs when
\begin{equation*}
\lim_{n\rightarrow\infty} \,[\dist(\alpha_n,\mathbb{S}_n)]^{1/n}\,=\,0\,.
\end{equation*}
Having in mind the method above, it is not hard to see that in this case, we have that
\begin{equation}\label{limitcase}
\limsup_{n\rightarrow\infty}\,e^{-\Re x}\,|L^{(\alpha_n)}_n(nx)|^{1/n}\,\leq \,e^{-s}\,,\;\,x\in \Gamma_s \setminus \{x_0(s)\}\,,
\end{equation}
\noindent for any $s>0\,.$ Thus, applying \cite[Theorem 5]{GaKu97}, \eqref{limitcase} implies that $\dsty \supp \mu_{\infty}\,\subset \,\overline{G_s}\,,$ for any $s>0\,.$ Since $\dsty \bigcap_{s>0}\,\overline{G_s}\,=\{0\}\,,$ the conclusion easily follows.

\section*{Acknowledgements}

R.O. thanks Professors A. B. J.
Kuijlaars, A. Mart\'{\i}nez Finkelshtein and H. Stahl for useful
discussions.

\end{document}